\documentclass[12pt]{article}
\usepackage{amsmath,amsfonts,epsfig}
\title{\bf On the Rank of a Coxeter Group}
\author{Michael L. Mihalik and John G. Ratcliffe\\ 
Mathematics Department, Vanderbilt University, \\
Nashville TN 37240, USA}
\newtheorem{theorem}{Theorem}[section]

\newtheorem{lemma}[theorem]{Lemma}

\newenvironment{proof}{{\bf Proof:\ }}{\hfill$\square$\vspace{.2in}}

\def\ov{\overline}

\date{}
\begin{document}
\maketitle

\vspace{.1in}


\section{Introduction} 
Let $W$ be a Coxeter group with Coxeter generators $S$. 
The {\it rank} of the Coxeter system $(W,S)$ is the cardinality $|S|$ of $S$. 
The Coxeter system $(W,S)$ has finite rank if and only if $W$ is finitely generated 
by Theorem 2(iii), Ch. IV, \S 1 of \cite{Bourbaki}. 
If $(W,S)$ has infinite rank, then $|S| = |W|$, 
since every element of $W$ is represented by a finite product of elements of $S$.  
Thus if $W$ is not finitely generated, the rank of $(W,S)$ is uniquely determined by $W$. 
If $W$ is finitely generated, then $W$ may have sets of Coxeter generators $S$ and $S'$ 
of different ranks. 
In this paper, we determine the set of all 
possible ranks for an arbitrary finitely generated Coxeter group $W$. 

This paper is a continuation of our previous paper with Steven Tschantz \cite{M-R-T} 
in which we studied the relationship between two sets $S$ and $S'$ of Coxeter generators 
of a finitely generated Coxeter group $W$. 
A {\it basic subset} of $S$ is 
a maximal subset $B$ of $S$ such that $B$ generates 
an irreducible, noncyclic, finite subgroup of $W$. 
In \cite{M-R-T}, we proved the Basic Matching Theorem 
which says that there is a natural bijection (matching) between 
the basic subsets of $S$ and the basic subsets of $S'$. 
A basic subset $B$ of $S$ matches a basic subset $B'$ of $S'$ 
if and only if $[\langle B\rangle,\langle B\rangle]$ is conjugate to 
$[\langle B'\rangle,\langle B'\rangle]$ in $W$. 
Usually matching basic subsets generate isomorphic groups, 
in which case, we say that the {\it basic subsets match isomorphically}\,; 
however, there are exceptions, due to well known 
isomorphisms between irreducible and reducible 
finite Coxeter groups (for instance the dihedral group 
${\bf D}_2(6)$ of order 12 and ${\bf A}_1\times {\bf A}_2$). 
We showed that nonisomorphic matching of basic subsets can be understood 
by {\it blowing up} Coxeter generating sets. 
This is a procedure to replace a given Coxeter generating set $S$ by a Coxeter 
generating set $R$ such that $|R| = |S|+1$. 
In \cite{M-R-T}, we proved that there exists a set of Coxeter generators $S'$ of $W$ 
such that a basic subset $B$ of $S$ matches a basic subset $B'$ of $S'$ 
with $|\langle B\rangle| > |\langle B'\rangle|$ if and only if $S$ can be blown up.  
We proved that $S$ has maximum rank over all sets of Coxeter generators of $W$ 
if and only if $S$ can not be blown up. 

In this paper, we study the reverse procedure of {\it blowing down} Coxeter generating sets, 
which was introduced by Mihalik in \cite{Mihalik}. 
We first determine necessary and sufficient conditions on $(W,S)$ 
such that there exists a set of Coxeter generators $S'$ of $W$ 
such that a basic subset $B$ of $S$ matches a basic subset $B'$ of $S'$ 
with $|\langle B\rangle| < |\langle B'\rangle|$. 
We then determine necessary and sufficient conditions on $(W,S)$ such that 
$W$ has a set of Coxeter generators $S'$ such that $|S'| < |S|$. 
As an application, we determine the rank spectrum of $W$.

\section{Preliminaries} 

A {\it Coxeter matrix} is a symmetric matrix $M = (m(s,t))_{s,t\in S}$ 
with $m(s,t)$ either a positive integer or infinity and $m(s,t) =1$ 
if and only if $s=t$. A {\it Coxeter system} with Coxeter matrix $M = (m(s,t))_{s,t\in S}$ 
is a pair $(W,S)$ consisting of a group $W$ and a set of generators $S$ for $W$ 
such that $W$ has the presentation
$$W =\langle S \ |\ (st)^{m(s,t)}:\, s,t \in S\ \hbox{and}\ m(s,t)<\infty\rangle$$
If $(W,S)$ is a Coxeter system with Coxeter matrix $M = (m(s,t))_{s,t\in S}$, 
then the order of $st$ is $m(s,t)$ for each $s,t$ in $S$, 
and so a Coxeter system determines its Coxeter matrix; moreover, any Coxeter matrix 
$M = (m(s,t))_{s,t\in S}$ determines a Coxeter system $(W,S)$ where $W$ 
is defined by the above presentation. If $(W,S)$ is a Coxeter system, 
then $W$ is called a {\it Coxeter group} and $S$ is called a set of {\it Coxeter generators} 
for $W$, and the cardinality of $S$ is called the {\it rank} of $(W,S)$.

Let $(W,S)$ be a Coxeter system. 
The {\it Coxeter diagram} ({\it {\rm C}-diagram}) of $(W,S)$ is the labeled undirected graph 
$\Delta(W,S)$ with vertices $S$ and edges 
$$\{(s,t) : s, t \in S\ \hbox{and}\ m(s,t) > 2\}$$
such that an edge $(s,t)$ is labeled by $m(s,t)$. 
A Coxeter system $(W,S)$ is said to be {\it irreducible} 
if its C-diagram is connected. 

A {\it visible subgroup} of $(W,S)$ 
is a subgroup of $W$ of the form $\langle A\rangle$ for some $A \subseteq S$. 
A visible subgroup $\langle A\rangle$ of $(W,S)$ is said to be {\it irreducible} 
if $(\langle A\rangle, A)$ is irreducible. 
A subset $A$ of $S$ is said to be {\it irreducible} if $\langle A\rangle$ is irreducible. 
A subset $A$ of $S$ is said to be a {\it component} of $S$ if $A$ is a maximal irreducible 
subset of $S$ or equivalently if $\Delta(\langle A\rangle, A)$ is a connected component 
of $\Delta(W,S)$.

The {\it presentation diagram} ({\it {\rm P}-diagram}) of $(W,S)$ is the labeled undirected graph 
$\Gamma(W,S)$ with vertices $S$ and edges 
$$\{(s,t) : s, t \in S\ \hbox{and}\ m(s,t) < \infty\}$$
such that an edge $(s,t)$ is labeled by $m(s,t)$.

We continue with the terminology of \cite{M-R-T}. 
In particular, we use Coxeter's notation on p. 297 of \cite{Coxeter}
for the irreducible spherical simplex reflection groups except 
that we denote the dihedral group ${\bf D}_2^k$ by ${\bf D}_2(k)$. 
Subscripts denote the rank of a Coxeter system in Coxeter's notation. 
Coxeter's notation partly agrees with but differs from Bourbaki's notation on p. 193 of 
\cite{Bourbaki}. 
Coxeter proved that every finite irreducible Coxeter system 
is isomorphic to exactly one 
of the Coxeter systems ${\bf A}_n$, $n\geq 1$, ${\bf B}_n$, $n\geq 4$, ${\bf C}_n$, 
$n\geq 2$, ${\bf D}_2(k)$, $k\geq 5$, 
${\bf E}_6$, ${\bf E}_7$, ${\bf E}_8$, ${\bf F}_4$, ${\bf G}_3$, ${\bf G}_4$.  
See \S 3 of \cite{M-R-T} for definitions. 
For uniformity of notation, 
we define ${\bf B}_3 = {\bf A}_3$, ${\bf D}_2(3) = {\bf A}_2$ and ${\bf D}_2(4) = {\bf C}_2$. 

Let $(W,S)$ be a Coxeter system. 
A {\it basic subset} of $S$ is a maximal irreducible subset $B$ of $S$ 
such that $\langle B\rangle$ is a noncyclic finite group. 
If $B$ is a basic subset of $S$, we call $B$ a {\it base} of $(W,S)$ 
and $\langle B\rangle$ a {\it basic subgroup} of $W$.

\begin{theorem} {\rm (Basic Matching Theorem, Theorem 4.18 \cite{M-R-T})} 
Let $W$ be a finitely generated Coxeter group with 
two sets of Coxeter generators $S$ and $S'$. 
Let $B$ be a base of $(W,S)$. 
Then there is a unique irreducible subset $B'$ of $S'$ such that 
$[\langle B\rangle,\langle B\rangle]$ is conjugate to 
$[\langle B'\rangle,\langle B'\rangle]$ in $W$. Moreover, 

\begin{enumerate}
\item the set $B'$ is a base of $(W,S')$, and we say that $B$ and $B'$ match,  

\item if $|\langle B\rangle|=|\langle B'\rangle|$, then $B$ and $B'$ have the same type 
and there is an isomorphism $\phi:\langle B\rangle \to \langle B'\rangle$ 
that restricts to conjugation on $[\langle B\rangle,\langle B\rangle]$ 
by an element of $W$, and we say that $B$ and $B'$ match isomorphically, 

\item if $|\langle B\rangle|<|\langle B'\rangle|$, then either
$B$ has type ${\bf B}_{2q+1}$ and 
$B'$ has type ${\bf C}_{2q+1}$ for some $q\geq 1$ or 
$B$ has type ${\bf D}_2(2q+1)$ and 
$B'$ has type ${\bf D}_2(4q+2)$ for some $q\geq 1$. 
Moreover, there is a monomorphism $\phi:\langle B\rangle \to \langle B'\rangle$ 
that restricts to conjugation on $[\langle B\rangle,\langle B\rangle]$ 
by an element of $W$. 
\end{enumerate}
\end{theorem}

\section{Blowing Down Coxeter Systems} 

Let $(W,S)$ be a Coxeter system of finite rank. 
In this section, we determine necessary and sufficient conditions 
for a base $B$ of $(W,S)$ to 
match a base $B'$ of $(W,S')$ with $|\langle B\rangle| < |\langle B'\rangle|$. 
If a base $B$ of $(W,S)$ matches a base $B'$ of $(W,S')$ 
with $|\langle B\rangle| < |\langle B'\rangle|$, 
then either $B$ is of type ${\bf B}_{2q+1}$ and $B'$ is of type ${\bf C}_{2q+1}$ 
for some $q\geq 1$ 
or $B$ is of type ${\bf D}_2(2q+1)$ and $B'$ is of type ${\bf D}_2(4q+2)$ 
for some $q\geq 1$ by the Basic Matching Theorem. 

If $a\in S$, the {\it neighborhood} of $a$ in P-diagram of $(W,S)$ is 
defined to be the set 
$N(a) =\{s\in S: m(s,a) < \infty\}$.  
If $A\subseteq S$, define 
$$A^\perp=\{s\in S: m(s,a) = 2\ \hbox{for all}\ a\in A\}.$$ 

The following lemma generalizes Proposition 2 of \cite{Mihalik}. 

\begin{lemma}\label{First Lemma} 
Let $B =\{x,y\}$ be a base of $(W,S)$ of type ${\bf D}_2(2q+1)$ 
that matches a base $B'$ of $(W,S')$ of type ${\bf D}_2(4q+2)$ for some $q\geq 1$. 
Then $N(x)\cap N(y) = B\cup B^\perp$. 
\end{lemma}
\begin{proof}
Suppose $s\in S-B$ with $m(s,x), m(s,y) <\infty$. 
Let $M\subseteq S$ be a maximal simplex containing $\{s,x,y\}$. 
Then there is a unique maximal simplex $M'\subseteq S'$ such 
that $\langle M\rangle$ is conjugate to $\langle M'\rangle$ by Prop. 4.21 of \cite{M-R-T}. 
By conjugating $S'$, we may assume that $\langle M\rangle=\langle M'\rangle$. 
Then $M'$ contains $B'$ and $[\langle B\rangle,\langle B\rangle]$ 
is conjugate to $[\langle B'\rangle,\langle B'\rangle]$ in $\langle M'\rangle$ 
by the Basic Matching Theorem. 

Let $B'=\{a,b\}$. Then $m(s',a)=m(s',b)=2$ for all $s'\in M'-B'$ by Theorem 8.7 of \cite{M-R-T}. 
Hence $B'$ is a component of $M'$. 
Therefore $[\langle B'\rangle,\langle B'\rangle]$ is a normal subgroup of $\langle M'\rangle$. 
Hence $[\langle B\rangle,\langle B\rangle]$ is a normal subgroup of $\langle M\rangle$. 
Therefore $s\{x,y\}s =\{x,y\}$ by Lemma 4.17 of \cite{M-R-T},   
and so $sxs=x$ and $sys=y$ by the deletion condition. 
Hence $N(x)\cap N(y) = B\cup B^\perp$. 
\end{proof}

The C-diagram of ${\bf B}_n$, for $n\geq 5$, is a Y-shaped diagram 
with $n$ vertices $b_1,\ldots, b_n$   
and two short arms.  We call the endpoints $b_{n-1}$ and $b_n$ of 
the short arms the {\it split ends} of the C-diagram of ${\bf B}_n$. 
The {\it split ends} of ${\bf B}_3$ are the end points 
$b_2$ and $b_3$ of the C-diagram of ${\bf B}_3 = {\bf A}_3$. 
The C-diagram of ${\bf C}_n$ is a linear diagram with $n$ vertices 
$c_1,\ldots, c_n$ and  
$m(c_i,c_{i+1})=3$ for $i=1,\ldots, n-2$ and $m(c_{n-1},c_n)=4$. 
Note that $b_i=c_i$ for $i=1,\ldots, n-1$, and $b_n=c_nc_{n-1}c_n$, 
and $b_{n-1}b_n = (c_{n-1}c_n)^2$. See \S 3 of \cite{M-R-T} for details. 

\begin{lemma}\label{Second Lemma} 
Let $\phi:{\bf B}_n \to {\bf C}_n$ be a monomorphism with $n$ odd and $n\geq 3$. 
Then $\phi$ maps $b_{n-1}b_n$ to a conjugate of $(c_{n-1}c_n)^2$ in ${\bf C}_n$. 
\end{lemma}
\begin{proof} 
Now $\phi({\bf B}_n)$ does not contain the center of ${\bf C}_n$, 
since $Z({\bf B}_n) = \{1\}$. 
Therefore either $\phi({\bf B}_n)={\bf B}_n$ or $\phi({\bf B}_n)=\theta({\bf B}_n)$ 
where $\theta$ is the automorphism of ${\bf C}_n$ defined by 
$\theta(c_i) = -c_i$, for $i=1,\ldots,n-1$ and $\theta(c_n)=c_n$. 
Now $\theta$ restricts to the identity on $[{\bf C}_n,{\bf C}_n]$, 
and so by composing $\phi$ with $\theta$ in the latter case, 
we may assume that $\phi({\bf B}_n)={\bf B}_n$. 
Every automorphism of ${\bf B}_n$ is inner by Theorem 31 of \cite{F-H}.  
Hence $\phi$ restricts to conjugation on $[{\bf B}_n,{\bf B}_n]$ 
by an element of ${\bf B}_n$. 
As $b_{n-1}b_n$ is in $[{\bf B}_n,{\bf B}_n]$, and $b_{n-1}b_n=(c_{n-1}c_n)^2$, 
we conclude that $\phi(b_{n-1}b_n)$ is conjugate to $(c_{n-1}c_n)^2$ in ${\bf C}_n$.
\end{proof}

Let $W$ be a finitely generated Coxeter group with two sets $S$ and $S'$   
of Coxeter generators, and let $A$ be a subset of $S$. 
Let $\ov A$ be the intersection of all subsets $B$ of $S$ such that $B$ contains $A$ 
and $\langle B\rangle$ is conjugate to $\langle B'\rangle$ for some $B'\subseteq S'$. 
Then $\ov A$ is the smallest subset $B$ of $S$ such that $B$ contains $A$ 
and $\langle B\rangle$ is conjugate to $\langle B'\rangle$ for some $B'\subseteq S'$ 
by Prop. 4.14 of \cite{M-R-T}. 
If $A$ is a spherical simplex, then $\ov A$ is a spherical simplex,  
since for any maximal spherical simplex $M\subseteq S$ that contains $A$,  
there exists a unique maximal spherical simplex $M'\subseteq S'$ such that 
$\langle M\rangle$ is conjugate to $\langle M'\rangle$ by Prop. 4.13 of \cite{M-R-T}.

\begin{lemma}\label{Third Lemma} 
Let $B$ be a base of $(W,S)$ of type ${\bf B}_{2q+1}$ 
that matches a base $B'$ of $(W,S')$ of type ${\bf C}_{2q+1}$ for some $q\geq 1$. 
Let $x,y$ be the split ends of the C-diagram of $(\langle B\rangle, B)$. 
Then $\ov{\{x,y\}} =\ov B$ and $N(x)\cap N(y) = B\cup B^\perp$. 
\end{lemma}
\begin{proof}
Let $C = \ov{\{x,y\}}$. 
Then $C$ is a spherical simplex of $(W,S)$ and $\langle C\rangle$ 
is conjugate to $\langle C'\rangle$ for some $C'\subseteq S'$. 
By conjugating $S'$, we may assume that $\langle C\rangle = \langle C'\rangle$. 
Let $a, b, c$ be the elements of $B'$ such that $m(a,b)=4$ and $m(b,c)=3$. 
Now $xy\in [\langle B\rangle,\langle B\rangle]$, 
and so $xy$ is conjugate to $(ab)^2$ by the 
Basic Matching Theorem and Lemma \ref{Second Lemma}. 
Hence there is a $w\in W$ such that $w(ab)^2w^{-1}\in \langle C'\rangle$. 
Now $\langle (ab)^2\rangle =[\langle a,b\rangle,\langle a,b\rangle]$. 
Let $u$ be the shortest element of $\langle C'\rangle w\langle a,b\rangle$. 
Then $u\{a,b\}u^{-1}\subseteq C'$ by Lemma 4.17 of \cite{M-R-T}. 
As $m(a,b)=4$, we deduce that $\{a,b\}\subseteq C'$ by Lemma 4.9 of \cite{M-R-T}. 
Hence $B'\subseteq C'$ by Lemma 8.1 of \cite{M-R-T}. 
By the Basic Matching Theorem, $B\subseteq C$. 
Hence $\ov B\subseteq C$.  As $\{x,y\}\subset B$, we have $C\subseteq \ov B$. 
Thus $C = \ov B$.

Suppose $s\in S-B$ with $m(s,x), m(s,y)<\infty$. 
Let $M\subseteq S$ be a maximal simplex containing $\{s,x,y\}$. 
Then there is a maximal simplex $M'\subseteq S'$ such 
that $\langle M\rangle$ is conjugate to $\langle M'\rangle$ by Prop. 4.21 of \cite{M-R-T}. 
By conjugating $S'$, we may assume that $\langle M\rangle=\langle M'\rangle$. 
Now $\ov{\{x,y\}}\subseteq M$, and so $B\subseteq M$. 
Hence $B'\subseteq M'$ by the Basic Matching Theorem. 
Moreover $m(s',t') = 2$ 
for all $(s',t')\in (M'-B')\times B'$ by Theorem 8.2 of \cite{M-R-T}. 
Hence $B'$ is a component of $M'$. 
Therefore $[\langle B'\rangle,\langle B'\rangle]$ is a normal subgroup of $\langle M'\rangle$. 
By the Basic Matching Theorem, $[\langle B\rangle,\langle B\rangle]$ 
is conjugate to $[\langle B'\rangle,\langle B'\rangle]$ in $\langle M\rangle$. 
Therefore $[\langle B\rangle,\langle B\rangle]$ is a normal subgroup of $\langle M\rangle$. 
Then $sBs = B$ by Lemma 4.17 of \cite{M-R-T}, 
and so $sts=t$ for all $t\in B$ by the deletion condition. 
Hence $N(x)\cap N(y) = B\cup B^\perp$. 
\end{proof}

The {\it odd diagram} of $W$ is the labeled undirected diagram $\Omega(W,S)$ 
obtained from the P-diagram of $(W,S)$ be deleting the even labeled edges. 
If $a\in S$, we define ${\rm Odd}(a)$ to be the vertex set of 
the connected component of $\Omega(W,S)$ containing $a$. 
By Prop. 3, Ch. IV, \S 1 of Bourbaki \cite{Bourbaki}, we have that 
$${\rm Odd}(a) = \{s\in S: s\ \hbox{is conjugate to}\ a\ \hbox{in}\ W\}.$$
We define 
$${\rm EOdd}(a) = {\rm Odd}(a)\cup \{s\in S: m(s,b)\ \hbox{is even for some}\ b\in {\rm Odd}(a)\}.$$

The next lemma has its genesis in Proposition 3 of \cite{Mihalik}.

\begin{lemma}\label{Fourth Lemma} 
Let $B$ be a base of $(W,S)$ that matches a base $B'$ of $(W,S')$ 
with $|\langle B\rangle| < |\langle B'\rangle|$. 
Then there exists $r\in \ov{B} - B$ such that $N(r) =B\cup B^\perp$  
and ${\rm Odd}(r) = \{r\}$, and if $K$ is the component of $B^\perp$ containing $r$, 
then $K$ is of type ${\bf A}_1$, ${\bf C}_{2q+1}$, or ${\bf D}_2(4q+2)$ for some $q\geq 1$; 
moreover, if $K\neq\{r\}$, then $K$ is a basic subset of $S$ and 
if $K'$ is the basic subset of $S'$ that matches $K$, 
then $K'$ is a component of $(B')^\perp$ and $K'\cup (K')^\perp = B'\cup(B')^\perp$.  
\end{lemma}
\begin{proof}
Let $C=\ov B$. 
Then $C$ is a spherical simplex of $(W,S)$ 
and $\langle C\rangle$ is conjugate to $\langle C'\rangle$ for some $C'\subseteq S'$. 
By conjugating $S'$, we may assume that $\langle C\rangle = \langle C'\rangle$. 
Then $C'$ contains $B'$ by the Basic Matching Theorem. 
Hence $B$ is a proper subset of $C$, since otherwise 
$\langle B'\rangle \subseteq \langle C'\rangle =\langle C\rangle = \langle B\rangle$ 
which is not the case, since $|\langle B\rangle|< |\langle B'\rangle|$.

Let $M\subseteq S$ be a maximal spherical simplex containing $B$, 
and let $M'\subseteq S'$ be the maximal spherical simplex 
such that $\langle M\rangle$ is conjugate to $\langle M'\rangle$. 
Then $M'$ contains $B'$ by the Basic Matching Theorem. 
Let $w$ be an element of $W$ such that $w\langle M\rangle w^{-1} = \langle M'\rangle$. 
Now $B'\subseteq \langle C\rangle \subseteq \langle M\rangle$. 
Hence $wB'w^{-1} \subset \langle M'\rangle$. 
Let $u$ be the shortest element of $\langle M'\rangle w\langle B'\rangle$. 
Then $uB'u^{-1}\subseteq M'$ by Lemma 4.3 of \cite{M-R-T}. 
As $B'$ is a base of $(W,S')$, we have that $uB'u^{-1} = B'$ by Lemma 4.10 of \cite{M-R-T},  
and so $u$ acts as a graph automorphism on $\langle B'\rangle$. 
Let $z'$ be the longest element of $\langle B'\rangle$. 
Then $uz'u^{-1} = z'$. 
Now $w = xuy$ with $x\in \langle M'\rangle$ and $y\in \langle B'\rangle$. 
Hence $wz'w^{-1} = z'$, since $z'$ is in the center of $\langle M'\rangle$. 
Therefore $z'$ is in the center of $\langle M\rangle = w^{-1}\langle M'\rangle w$. 
As $B\cup B^\perp$ is the union of all the maximal spherical simplices of $(W,S)$ 
that contain $B$, we deduce the $z'$ is in the center of $B^\perp$. 
Hence, there are distinct components $K_1,\ldots, K_n$ of $B^\perp$,  
with nontrivial center, such that $z' = z_1\cdots z_n$ 
with $z_i$ the longest element of $\langle K_i\rangle$ for each $i=1,\ldots,n$.  
As $z'\in \langle C\rangle$, we have that $K_i\subseteq C$ for each $i=1,\ldots, n$ 
by Prop. 7, Ch. IV, \S 1 of \cite{Bourbaki}, since every reduced form of $z_i$ 
involves every element of $K_i$ for each $i=1,\ldots,n$. 

Define a homomorphism $\rho:\langle C'\rangle \to \langle z'\rangle$ as follows. 
To begin with, define $\rho(s') = 1$ if $s'\in C'-B'$.  
By the Basic Matching Theorem, $B'$ is of type ${\bf C}_{2p+1}$ or ${\bf D}_2(4p+2)$ 
for some $p\geq 1$. 
If $B'$ is of type ${\bf C}_{2p+1}$, let $a'\in B'$ be such that $B'\cap {\rm Odd}(a') =\{a'\}$, 
and define $\rho(a') = z'$ and $\rho(s') = 1$ for each $s'\in B'-\{a'\}$.  
Suppose $B'$ is of type ${\bf D}_2(4p+2)$ and $B' = \{a',b'\}$. 
By Lemma 8.6 of \cite{M-R-T}, one of $a'$ or $b'$, say $a'$, has the property 
that if $a'\in A'\subseteq S'$ and $\langle A'\rangle$ is conjugate to $\langle A\rangle$ 
for some $A\subseteq S$, then $B'\subseteq A'$. 
Define $\rho(a') = z'$ and $\rho(b') =1$. 
In both cases, $\rho$ is well defined and $\rho(z') = z'$. 

As $z'=z_1\cdots z_n$, there is an $i$ such that $\rho(z_i) = z'$. 
By reindexing, we may assume $i=1$. 
Let $K = K_1$. 
Then there exists $r\in K$ such that $\rho(r) = z'$. 
As $r\in \ov B$, we have that $\ov{\{r\}}\subseteq \ov B$. 
Let $A= \ov{\{r\}}$. 
Then $\langle A\rangle$ is conjugate in $\langle C'\rangle$
to $\langle A'\rangle$ for some $A' \subseteq C'$ by Prop. 4.14 of \cite{M-R-T}. 
Now $\rho(\langle A'\rangle) = \rho(\langle A\rangle) = \langle z'\rangle$. 
Hence $a' \in A'$. 
Then $B'\subseteq A'$ by Lemma 8.1 of \cite{M-R-T} or the choice of $a'$.  
Hence $B \subseteq A$ by the Basic Matching Theorem,  
and so $\ov B\subseteq A$. 
Hence $\ov{\{r\}}= \ov B$. 

As $K$ is a component of $B^\perp$, we have $B\cup B^\perp \subseteq N(r)$. 
Suppose $s\in N(r)$. 
Let $M\subseteq S$ be a maximal spherical simplex containing $\{r,s\}$. 
Then there is a maximal spherical simplex $M'\subseteq S'$ 
such that $\langle M\rangle$ is conjugate to $\langle M'\rangle$. 
Then $\ov{\{r\}}\subseteq M$, and so $B\subseteq M$. 
Therefore $s\in B\cup B^\perp$, since $B$ is a basic subset of $M$.  
Thus $N(r) = B\cup B^\perp$. 

As $[\langle K\rangle, \langle K\rangle]$ is in the kernel of $\rho$, 
we have that $z_1$ is not in $[\langle K\rangle, \langle K\rangle]$. 
Therefore $K$ is of type ${\bf A}_1$, ${\bf C}_{2q+1}$, ${\bf D}_2(4q+2)$, ${\bf E}_7$, 
or ${\bf G}_3$ for some $q\geq 1$. 
Suppose $K$ is of type ${\bf C}_{2q+1}$ for some $q\geq 1$. 
Let $a\in K$ be such that $K\cap{\rm Odd}(a) = \{a\}$. 
Then $a[\langle K\rangle,\langle K\rangle] = z_1[\langle K\rangle,\langle K\rangle]$. 
As the restriction of $\rho$ to $\langle K\rangle$ factors through 
$\langle K\rangle/[\langle K\rangle,\langle K\rangle]$, we may assume that $r=a$. 
Then ${\rm Odd}(r) = \{r\}$. 

If $K= \{r\}$, then we are done. 
Suppose $K\neq \{r\}$. 
Then $K$ is a basic subset of $S$, since $N(r) = B\cup B^\perp$ 
and $K$ is a component of $B^\perp$. 
Let $K'$ be the basic subset of $C'$ that matches $K$. 
Then $K'$ is the basic subset of $S'$ that matches $K$ 
by the Basic Matching Theorem. 
Let $M'\subseteq S'$ be a maximal spherical simplex that contains $B'$, 
and let $M\subseteq S$ be the maximal spherical simplex such that $\langle M\rangle$ 
is conjugate to $\langle M'\rangle$. 
Then $M$ contains $B$ by the Basic Matching Theorem. 
Now $K\subseteq C\subseteq M$ and $K$ is a basic subset of $M$. 
Therefore $K'$ is a basic subset of $M'$ by the Basic Matching Theorem. 
Hence $K'$ is a component of $M'$.  
Therefore $K'$ is a component of $(B')^\perp$. 

Suppose $s'\in (K')^\perp$. 
Let $M'\subseteq S'$ be a maximal spherical simplex that contains $K'\cup\{s'\}$, 
and let $M\subseteq S$ be the maximal spherical simplex such that $\langle M\rangle$ 
is conjugate to $\langle M'\rangle$. 
Then $K\subseteq M$ by the Basic Matching Theorem. 
Hence $M$ contains $r$, and so $B\subseteq \ov{\{r\}}\subseteq M$. 
Therefore $M'$ contains $B'$ by the Basic Matching Theorem. 
Hence $s'\in B'\cup (B')^\perp$. 
Therefore $(K')\cup (K')^\perp = B'\cup(B')^\perp$. 
If $K$ is of type ${\bf C}_{2q+1}$ or ${\bf D}_2(4q+2)$, we are done.

Suppose $K$ is of type ${\bf E}_7$ or ${\bf G}_3$. 
Then $\rho(s) = z'$ for each $s$ in $K$, since $s$ and $r$ are conjugate in $\langle K\rangle$ 
for each $s\in K$. 
Hence $\ov{\{s\}} = \ov B$ and $N(s) =B\cup B^\perp$ for each $s\in K$ by the above argument.  
Therefore ${\rm Odd}(r) = K$ and ${\rm EOdd}(r) = B\cup B^\perp$. 
Suppose $n > 1$. 
By Lemma 28 of \cite{H-M}, there is an automorphism $\theta$ of $W$ 
such that $\theta(s) = s$ for all $s\in S-K$, 
and $\theta(s) = sz_2\cdots z_n$ for all $s\in K$. 
The longest element of $\langle\theta(K)\rangle$ is $z_1\cdots z_n=z'$. 
Now replace $S$ by $\theta(S)$. 
Let $\ov{\theta(B)}$ denote $\ov B$ with respect to the Coxeter generators $\theta(S)$ and $S'$. 
Then $\ov{\theta(B)}\subseteq \theta(C)$, 
since $\langle \theta(C)\rangle = \langle C\rangle =\langle C'\rangle$. 
If $\ov{\theta(B)}$ is a proper subset of $\theta(C)$, we return to the start of the proof. 
As $C$ is finite, we will eventually be done 
or have $\ov{\theta(B)}=\theta(C)$. 
Thus we may assume without loss of generality that $n = 1$ and $z' = z_1$. 

Let $\ell'$ be the longest element of $\langle K'\rangle$. 
Define a homomorphism $\eta:\langle C'\rangle \to \langle \ell'\rangle$ as follows. 
Define $\eta(s') = 1$ for all $s'\in C'-K'$ 
and define $\eta(s') = \ell'$ for all $s' \in K'$. 
Then $\eta$ is well defined and $\eta(\ell') = \ell'$.  
By the same argument as above,  
$\ell'$ is in the center of $B^\perp$ 
and there is a component $L$ of $B^\perp$ 
such that $L \subseteq C$ and $L$ has nontrivial center, 
and if $\ell$ is the longest element of $\langle L\rangle$, then $\eta(\ell) = \ell'$ 
and $\eta(t) = \ell'$ for some $t\in L$. 
As $z'\in B'\subseteq C'-K'$, we have that $\eta(z') = 1$, and so $L\neq K$.

Let $A = \ov{\{t\}}$. 
Then $\langle A\rangle$ is conjugate in $\langle C'\rangle$
to $\langle A'\rangle$ for some $A' \subseteq C'$ by Prop. 4.14 of \cite{M-R-T}. 
Now $\eta(\langle A'\rangle) = \eta(\langle A\rangle) = \langle \ell'\rangle$. 
Hence $s'\in A'$ for some $s'\in K'$. 
Now killing $s'$ in $\langle C'\rangle$ kills $K'$.  
As $[\langle K\rangle, \langle K\rangle] = [\langle K'\rangle,\langle K'\rangle]$  
and $\langle A\rangle$ is conjugate to $\langle A'\rangle$ in $\langle C\rangle$,  
the group $\langle A\rangle$ contains an element that kills 
$[\langle K\rangle, \langle K\rangle]$ in $\langle C\rangle$. 
Therefore $A$ contains an element $s\in K$. 
Now $\ov{\{s\}} \subseteq A$. 
As $\ov{\{s\}} = \ov B$, we have $B\subseteq A$. 
Therefore $\ov{\{t\}} = \ov B$. 
As before, $N(t) = B\cup B^\perp$. 

By the same argument as above, $L$ is of type 
${\bf A}_1$, ${\bf C}_{2q+1}$, ${\bf D}_2(4q+2)$, ${\bf E}_7$, or ${\bf G}_3$, 
and if $L$ is of type ${\bf A}_1$, ${\bf C}_{2q+1}$, or ${\bf D}_2(4q+2)$, 
we are done.  
Suppose $L$ is of type  ${\bf E}_7$ or ${\bf G}_3$.  
Then by the same argument as above, 
${\rm Odd}(t) = L$ and ${\rm EOdd}(t) = B\cup B^\perp$,  
and we may assume that $\ell = \ell'$. 
By Lemma 38 of \cite{H-M}, there is an automorphism $\beta$ of $W$ 
such that $\beta(s) = s$ for all $s\in S-(K\cup L)$, 
and $\beta(s) = s\ell z'$ for each $s \in K\cup L$. 
Then $\beta(\ell) = z'$ and $\beta(z') = \ell$. 
As $\beta$ fixes each element of $S-(K\cup L)$ and 
$\beta$ leaves $\langle C\rangle$ invariant, 
we may replace $S$ by $\beta(S)$.  
Then $K$ is replaced by $\beta(L)$ 
and $\beta(K)$ is removed as a possibility for replacing $L$, 
since $\langle \beta(K)\rangle = \langle K'\rangle$.  
In the above procedure only Coxeter generators of components of $B^\perp$ 
of type ${\bf E}_7$ or ${\bf G}_3$ are replaced. 
By repeating this procedure a finite number of times, 
we can remove the possibility that $L$ is of type ${\bf E}_7$ or ${\bf G}_3$, 
and we are done. 
\end{proof}

A {\it cycle} of $S$ is a sequence $\{c_1,\ldots, c_n\}$ of distinct elements of $S$ 
so that $m(c_i,c_{i+1}) <\infty$ for $i=1,\ldots,n-1$ and $m(c_n,c_1)<\infty$. 
A {\it chord} of cycle $C =\{c_1,\ldots, c_n\}$ of $S$ 
is a pair of distinct elements $c_i,c_j$ of $C$ 
such that $m(c_i,c_j) < \infty$ and $c_i, c_j$ 
are neither consecutive terms of the cycle nor the end terms of the cycle. 

The next lemma generalizes Proposition 4 of \cite{Mihalik}. 

\begin{lemma}\label{Fifth Lemma}  
Let $B$ be a base of $(W,S)$ that matches a base $B'$ of $(W,S')$  
with $|\langle B\rangle|<|\langle B'\rangle|$. 
If $B$ is of type ${\bf D}_2(2q+1)$, let $B =\{x,y\}$.  
If $B$ is of type ${\bf B}_{2q+1}$, let $\{x,y\}$ be the set of split ends of the C-diagram of 
$(\langle B\rangle, B)$.  
Then $\{x,y\}$ is not part of a chord-free cycle of $S$ of length at least 4. 
\end{lemma}
\begin{proof}
On the contrary, suppose $C\subseteq S$ is a chord-free cycle of length at least 4 
that contains $\{x,y\}$.  
We may assume that $|S|$ is as small as possible. 
By Lemmas \ref{First Lemma} and \ref{Third Lemma}, we have that
 $\ov{\{x,y\}} = \ov B$ and $N(x)\cap N(y) = B\cup B^\perp$. 
By Lemma \ref{Fourth Lemma}, 
there is an $r\in \ov B-B$ such that $N(r)=B\cup B^\perp$ and ${\rm Odd}(r)=\{r\}$. 
Now $C\cap (B\cup B^\perp) = \{x,y\}$, since $C$ is cord-free of length at least 4. 
Let $a \in C-B$. 
Then $(B\cup B^\perp)-\{r\}$ is an $(a,r)$-separator of $S$, 
that is, every path from $a$ to $r$ in the P-diagram of $(W,S)$ passes 
through $(B\cup B^\perp)-\{r\}$.

Let $S_0$ be a c-minimal separator (see \S 6 of \cite{M-R-T}) 
of $S$ such that $S_0$ is conjugate to a subset 
of $(B\cup B^\perp)-\{r\}$.  
By Lemma 4.9 of \cite{M-R-T}, we have that
$S_0=S_1\cup S_2$ with $S_1$ a spherical simplex, 
$S_2\subseteq S_1^\perp$, 
and $wS_2w^{-1} \subseteq S$ if and only if $w=1$. 
Then $S_2\subseteq (B\cup B^\perp)-\{r\}$. 
By Theorem 6.1 of \cite{M-R-T}, there exists $S_0'\subseteq S'$, 
a reduced visual graph of groups decomposition $\Lambda$ for $(W,S)$, 
and a reduced visual graph of groups decomposition $\Lambda'$ for $(W,S')$ 
such that $\langle S_0\rangle$ is conjugate to $\langle S_0'\rangle$, 
and the edge groups of $\Lambda$ and $\Lambda'$ are conjugate to $\langle S_0\rangle$, 
and there is a 1-1 correspondence between the vertices of $\Lambda$ and 
the vertices of $\Lambda'$ such that each vertex group of $\Lambda$ is conjugate 
to the corresponding vertex group of $\Lambda'$. 

Now $r$ is not in an edge group of $\Lambda$, 
since $r$ is not conjugate to an element of $(B\cup B^\perp)-\{r\}$. 
Let $V$ be the vertex group of $\Lambda$ that contains $r$. 
Then $N(r)\subset V$, and so $B\cup B^\perp \subset V$. 
Now $\{x,y\}$ is not contained in an edge group $E$ of $\Lambda$, since otherwise 
$r\in \ov{\{x,y\}} \subset E$. 
We claim that $C\subset V$. 
On the contrary, suppose $C\not\subset V$.  
Let $C=\{c_1,\ldots,c_n\}$ with $x = c_1$ and $c_n=y$, and 
$m(c_i,c_{i+1})<\infty$ for each $i=1,\ldots,n-1$. 
Let $k$ be the first index such that $c_k\not\in V$ 
and let $\ell$ be the last index such that $c_\ell\not\in V$.  
Then $c_{k-1}$ is in an edge group $E$ of $\Lambda$ 
that is a subgroup of $V$ and $c_{\ell+1}$ is in an edge group $F$ of $\Lambda$ 
that is a subgroup of $V$.  
Now $E=F$, since the graph of $\Lambda$ is a tree. 
As $\{x,y\}\not\subset E$, we have that $\{c_{k-1},c_{\ell+1}\}\neq\{x,y\}$. 
Now $E$ is conjugate to $\langle S_0\rangle$, and so there exists $S_3\subseteq S_2^\perp$ 
such that $E = \langle S_3\cup S_2\rangle$ and $S_3$ is conjugate to $S_1$. 
Hence $S_3$ is a spherical simplex. 
Now $(C-B)\cap (B\cup B^\perp) = \emptyset$, since $C$ is chord-free. 
As $S_2\subseteq B\cup B^\perp$, we deduce that $\{c_{k-1},c_{\ell+1}\} \cap S_3 \neq \emptyset$, 
and so $c_{k-1}$ and $c_{\ell+1}$ are joined by a chord, which is a contradiction. 
Thus $C\subset V$. 

By conjugating $S'$, we may assume that $V$ is a vertex group $V'$ of $\Lambda'$. 
Then $B'\subset V'$ by the Basic Matching Theorem. 
Now ${\rm rank}(V) < |S|$, and so we have a contradiction to the minimality of $|S|$. 
\end{proof}

\begin{lemma} {\rm (Blow-Down Lemma)}\label{Sixth Lemma}  
Let $(W,S)$ be a Coxeter system of finite rank, and 
let $B$ be a base of $(W,S)$ of type ${\bf B}_{2q+1}$ or ${\bf D}_2(2q+1)$ for some $q\geq 1$. 
If $|B| = 2$, let $B =\{x,y\}$.  
If $B$ is of type ${\bf B}_{2q+1}$, let $\{x,y\}$ be the set of split ends of the C-diagram of 
$(\langle B\rangle, B)$.  
Let $r\in B^\perp$ such that $N(r) = B\cup B^\perp$ and $\{r\}$ is a component of $B^\perp$. 
Suppose $N(y) = B\cup B^\perp$. 
Let $\ell$ be the longest element of $\langle B\rangle$, let $a=r\ell$,  
let $S' = (S-\{r,y\})\cup \{a\}$, and let $B'=(B-\{y\})\cup \{a\}$.  
Then $S'$ is a set of Coxeter generators for $W$ such that 
\begin{enumerate}
\item the set $B'$ is a base of $(W,S')$ that matches $B$ 
with $|\langle B\rangle|<|\langle B'\rangle|$, 
\item $(B')^\perp = B^\perp-\{r\}$,  
\item the neighborhood of $a$ satisfies $N(a) = B'\cup (B')^\perp$, 
\item the basic subsets of $S$ and $S'$ are the same except for $B$ and $B'$. 
\end{enumerate}
\end{lemma}
\begin{proof}
Consider the Coxeter presentation
$$W = \langle S\ |\ (st)^{m(s,t)}: s,t \in S\ \hbox{and}\ m(s,t)<\infty\rangle.$$
Let $\ell$ be the longest element of $\langle B\rangle$. Then $\ell^2=1$. 
Regard $\ell$ as a reduced word in the elements of $B$. 
Add the generator $a$ and the relation $a=r\ell$ to the above presentation of $W$. 
Then we can add the relators $a^2$ and $(sa)^2$ and $(as)^2$ for all $s\in B^\perp-\{r\}$. 
As $\ell s\ell = s$ for all $s\in B-\{x,y\}$, 
we can add the relators $(sa)^2$ and $(as)^2$ for all $s\in B-\{x,y\}$. 

Next delete the generator $r$ and the relation $a = r\ell$ and replace $r$ 
by $a\ell$ in the remaining relators. 
The relator $r^2$ is replaced by $(a\ell)^2$. 
We delete the relators $(sa\ell)^2$ and $(a\ell s)^2$ for $s\in B\cup B^\perp-\{r,x,y\}$, 
since they are equivalent to $(a\ell)^2$. 

As $\ell x\ell = y$ and $\ell a \ell = a$, we have 
$xa\ell xa\ell = xa\ell x\ell a\ell\ell = xaya$, 
and so we can replace $(xa\ell)^2$ by $xaya$.  
Likewise $(a\ell x)^2$ can be replaced by $ayax$, 
and $(ya\ell)^2$ can be replaced by $yaxa$, 
and $(a\ell y)^2$ can be replaced by $axay$. 

Next we delete the generator $y$ and the relators $xaya$, $ayax$, $ayax$, $axay$,   
and replace $y$ by $axa$ in the remaining relators. 
We have eliminated all the relators originally involving $r$ except for $(a\ell)^2$. 

We delete the relators $(saxa)^2$ and $(axas)^2$ for all $s\in B\cup B^\perp-\{r,x,y\}$, 
since they are equivalent to $(axa)^2 = ax^2a$. 
The relator $y^2$ is replaced by $(axa)^2$ which we can delete. 

Assume first that $B$ is of type ${\bf B}_{2q+1}$. 
The relators $(xy)^2$ and $(yx)^2$ are replaced by $(xa)^4$ and $(ax)^4$.  
Let $t\in B$ be such that $m(t,y) = 3$. Then $m(t,x)=3$. 
Now $(ty)^3 = (taxa)^3 = (atxa)^3 = a(tx)^3a$, 
and so $(ty)^3$ can be deleted. 
Likewise $(yt)^3$ can be deleted. 
The relator $(a\ell)^2$ can be deleted, since it is redundant. 
Then we obtain a Coxeter presentation for $W$ with Coxeter generators $S'$ 
and $(\langle B'\rangle, B')$ of type ${\bf C}_{2q+1}$. 

Now assume that $B$ is of type ${\bf D}_2(2q+1)$. 
The relators $(xy)^{2q+1}$ and $(yx)^{2q+1}$ are replaced by $(xa)^{4q+2}$ 
and $(ax)^{4q+2}$. 
The relator $(a\ell)^2$ can be deleted, since it is redundant. 
Then we obtain a Coxeter presentation for $W$ with Coxeter generators $S'$ 
and $(\langle B'\rangle, B')$ of type ${\bf D}_2(4q+2)$. 

As  $a = r\ell$ and $\{r\}$ is a component of $B^\perp$, 
we have that $B^\perp-\{r\} \subseteq (B')^\perp$. 
Suppose $s\in (B')^\perp$. Then $s\in S-(B\cup\{r\})$. 
As $s$ commutes with $a=r\ell$, we have that $s\in N(r) = B\cup B^\perp$ 
by Lemma 8.3 of \cite{M-R-T}. Hence $s\in B^\perp-\{r\}$. 
Therefore $(B')^\perp = B^\perp-\{r\}$. 

Clearly, we have $B'\cup (B')^\perp \subseteq N(a)$. 
Suppose $s \in N(a)-B'$. 
Then $s\in N(r) = B\cup B^\perp$ by Lemma 8.3 of \cite{M-R-T}. 
Hence $s \in B^\perp-\{r\} = (B')^\perp$,  
and so $N(a) = B'\cup (B')^\perp$.  
Therefore $B'$ is a base of $(W,S')$ 
and $B'$ is the only base of $(W,S')$ that contains $a$. 
The base $B'$ matches $B$, 
since $[\langle B'\rangle, \langle B'\rangle]  = [\langle B\rangle,\langle B\rangle]$. 
As $N(y) = B\cup B^\perp$, we have that $B$ is the only base of $(W,S)$ 
that contains $y$. 
Therefore the basic subsets of $S$ and $S'$ are the same except for $B$ and $B'$.  
\end{proof}

Let $B$ and $r\in B^\perp$ be as in the Blow-Down Lemma. 
We call $r$ a {\it sink} for the base $B$. 
Let $B'$ and $S'$ be as in the Blow-Down Lemma.  
We say that $(W,S')$ is obtained by {\it blowing down} $(W,S)$ 
along the base $B$. 
We also say that $B'$ has been obtained by {\it blowing down} $B$. 

The next theorem has its genesis in Proposition 5 of \cite{Mihalik}. 

\begin{theorem}\label{Blow-Down Theorem}  
{\rm (Blow-Down Theorem)} 
Let $(W,S)$ be a Coxeter system of finite rank, and 
let $B$ be a base of $(W,S)$ of type ${\bf B}_{2p+1}$ or ${\bf D}_2(2p+1)$ for some $p\geq 1$. 
If $|B| = 2$, let $B =\{x,y\}$.  
If $B$ is of type ${\bf B}_{2p+1}$, let $\{x,y\}$ be the set of split ends of the C-diagram of 
$(\langle B\rangle, B)$.  
Then $W$ has a set of Coxeter generators $S'$ such that $B$ matches a base $B'$ of $(W,S')$ 
with $|\langle B\rangle|<|\langle B'\rangle|$ if and only if 
\begin{enumerate}
\item the neighborhoods of $x$ and $y$ satisfy $N(x)\cap N(y) = B\cup B^\perp$,
\item the set $\{x,y\}$ is not part of a chord-free cycle of $S$ of length at least 4, 
\item there exists $r\in B^\perp$ such that $N(r) = B\cup B^\perp$ and ${\rm Odd}(r) = \{r\}$, 
and if $K$ is the component of $B^\perp$ containing $r$, then $K$ is of type 
${\bf A}_1$, ${\bf C}_{2q+1}$, or ${\bf D}_2(4q+2)$ for some $q\geq 1$.
\end{enumerate}
\end{theorem}
\begin{proof}
Suppose $W$ has a set of Coxeter generators $S'$ 
such that $B$ matches a base $B'$ of $(W,S')$ 
with $|\langle B\rangle|<|\langle B'\rangle|$. 
Then condition (1) follows from Lemmas \ref{First Lemma} and \ref{Third Lemma},  
condition (2) follows from Lemma \ref{Fifth Lemma}, and  
condition (3) follows from Lemma \ref{Fourth Lemma}. 

Conversely, suppose conditions (1), (2), and (3) are satisfied. 
Let $S_0 = B\cup B^\perp$, and let $T = S-S_0$. 
Let $T_x$ be the set of all $t\in T$ such that there is a 
sequence $t_1,\ldots, t_n$ in $T$ such that $m(x,t_1) < \infty$, 
$m(t_i,t_{i+1})<\infty$ for each $i=1,\ldots, n-1$, and $t_n=t$.  
Define $T_y$ similarly. 
We claim that $T_x\cap T_y = \emptyset$. 
On the contrary suppose that $T_x\cap T_y \neq \emptyset$. 
Then there is a cycle $C$ of $S$ such that $C\cap S_0 = \{x,y\}$. 
Assume that $C$ is as short as possible. 
Then $C$ is chord-free. 
By condition (1), we deduce that $C$ has length at least 4, 
but this contradicts condition (2). 
Therefore $T_x\cap T_y = \emptyset$. 

Let $S_1 = S-T_y$ and $S_2 = S_0\cup T_y$. 
Then $S=S_1\cup S_2$ and $S_1\cap S_2 = S_0$, 
and $m(a,b) = \infty$ for all $a\in S_1-S_0$ and $b\in S_2-S_0$. 
Let $\ell$ be the longest element of $\langle B\rangle$. 
Then $\ell S_0\ell^{-1} = S_0$ and 
the triple $(S_1,\ell, S_2)$ determines an elementary twist (see \S 5 of \cite{M-R-T})  
of $(W,S)$ giving a new Coxeter generating set $S_\ast = S_1\cup \ell S_2 \ell^{-1}$. 
As $\ell y\ell^{-1} = x$, 
we have $\ell T_y\ell^{-1}\subseteq T_x$ with respect to $S_\ast$,  
and so by replacing $S$ by $S_\ast$, 
we may assume $T_y = \emptyset$. 
Then $N(y) = B\cup B^\perp$.  

If $K = \{r\}$, then we can blow down $B$.  
Hence $W$ has a set of Coxeter generators $S'$ such that 
$B$ matches a base $B'$ of $(W,S')$ 
with $|\langle B\rangle|<|\langle B'\rangle|$ by Lemma \ref{Sixth Lemma}. 
If $K\neq \{r\}$, we can blow up $S$ along $K$ by Theorems 8.4 and 8.8 of \cite{M-R-T}. 
This creates a sink for $B$, 
which allows us to blow down $B$.  
Therefore $W$ has a set of Coxeter generators $S'$ such that $B$ matches a base $B'$ of $(W,S')$ 
with $|\langle B\rangle|<|\langle B'\rangle|$ by Lemma \ref{Sixth Lemma}. 
\end{proof}

The proof of the Blown-Down Theorem indicates that we may have to blow up along 
one base in order to create a sink before we can blow down along another base.  
For example, the base ${\bf D}_2(3)$ 
of the Coxeter system 
${\bf C}_3\times {\bf D}_2(3)$ can be blown down only after the system is blown up 
along the base ${\bf C}_3$ to yield the system ${\bf B}_3\times {\bf A}_1\times{\bf D}_2(3)$. 
Then ${\bf A}_1$ is a sink for the base ${\bf D}_2(3)$, and so now we can  
blow down ${\bf D}_2(3)$ to obtain the system ${\bf B}_3\times{\bf D}_2(6)$. 
If we blow up a Coxeter system and then blow down the resulting Coxeter system,  
the initial and final systems have the same rank. 
For example, the initial system ${\bf C}_3\times {\bf D}_2(3)$ has 
the same rank as the final system ${\bf B}_3\times{\bf D}_2(6)$. 

\section{Contracting Coxeter Systems}  

In this section, we determine necessary and sufficient conditions on $(W,S)$ such that $W$ 
has a set of Coxeter generators $S'$ such that $|S'|< |S|$. 

\begin{theorem}
{\rm (Contracting Theorem)} 
Let $(W,S)$ be a Coxeter system of finite rank.  
Then $W$ has a set of Coxeter generators $S'$ 
such that $|S'| < |S|$ if and only if 
there is a base $B$ of $(W,S)$ of type ${\bf B}_{2p+1}$ or ${\bf D}_2(2p+1)$ 
for some $p\geq	1$ satisfying conditions (1), (2), (3) of the Blow-Down Theorem 
with $K = \{r\}$. 
\end{theorem}
\begin{proof}
Suppose there is a base $B$ of $(W,S)$ of type ${\bf B}_{2p+1}$ or ${\bf D}_2(2p+1)$ 
for some $p\geq	1$ satisfying conditions (1), (2), (3) of the Blow-Down Theorem 
with $K = \{r\}$. 
By twisting $(W,S)$ as in the proof of the Blow-Down Theorem, 
leaving $B\cup B^\perp$ invariant, 
we may assume that $N(y)=B\cup B^\perp$. 
Then $W$ has a set of Coxeter generators $S'$ such that $|S'| = |S|-1$  
by Lemma \ref{Sixth Lemma}.   

Conversely, suppose $W$ has a set of Coxeter generators $S'$ such that $|S'| < |S|$. 
We may assume that $S'$ has the maximum possible number of basic subsets that isomorphically 
match basic subsets of $S$. 
Now $S$ has a basic subset $B$ that nonisomorphically matches a basic subset $B'$ of $S'$ 
by the Simplex Matching Theorem (Theorem 7.7 of \cite{M-R-T}). 
Then $|\langle B\rangle| \neq |\langle B'\rangle|$ by the Basic Matching Theorem. 

Assume first that $|\langle B\rangle| < |\langle B'\rangle|$. 
Then $B$ satisfies conditions (1), (2), (3) of the Blow-Down Theorem. 
Let $r$ and $K$ be as in Lemma \ref{Fourth Lemma}. 
If $K =\{r\}$, we are done. 
Suppose that $K$ is of type ${\bf C}_{2q+1}$ or ${\bf D}_2(4q+2)$ for some $q\geq 1$. 
Then $K$ is a basic subset of $S$.  
Let $K'$ be the basic subset of $S'$ that matches $K$. 
We claim that $K$ isomorphically matches $K'$. 
On the contrary, suppose that $K$ nonisomorphically matches $K'$. 
By Lemma \ref{Fourth Lemma}, we have that $K'$ is a component of $(B')^\perp$ and 
$(K')\cup (K')^\perp = B'\cup(B')^\perp$. 
By Theorems 8.4-8.8 of \cite{M-R-T}, we can blow up $S'$ along $B'$,  
and after twisting as in the proof of the Blow-Down Theorem, 
leaving $(K')\cup (K')^\perp$ invariant, we can blow down $K'$ 
to obtain a set of Coxeter generators $S''$ such that $|S''| = |S'|$ and $S''$ 
has two more basic subsets than $S'$ isomorphically matching basic subsets of $S$, 
which contradicts the choice of $S'$. 
Thus $K$ isomorphically matches $K'$.

Let $z'$ be the longest element of $\langle B'\rangle$. 
As in the proof of Lemma \ref{Fourth Lemma}, the element $z'$ is in the center of $B^\perp$.  
By applying the automorphism $\theta$ of $W$ defined by $\theta(s) = s$ for all $s\in S-\{r\}$ 
and $\theta(r) = rz_2\cdots z_n$ as in the proof of Lemma \ref{Fourth Lemma}, if necessary, 
we may assume that $z'$ is the longest element of $\langle K\rangle$.     
By Theorems 8.5 and 8.7 of \cite{M-R-T}, the group $W$ has a set of Coxeter generators $S''$ 
such that $K$ matches a base $K''$ of $(W,S'')$ with $|\langle K\rangle| > |\langle K''\rangle|$. 
Therefore $W$ has a set of Coxeter generators $S''$ 
such that $K'$ matches a base $K''$ of $(W,S'')$ with $|\langle K'\rangle| > |\langle K''\rangle|$. 
Let $\ell'$ be the longest element of $\langle K'\rangle$. 
If $K'$ is of type ${\bf C}_{2q+1}$, let $a'\in K'$ be such that $K'\cap{\rm Odd}(a') = \{a'\}$. 
If $K'$ is of type ${\bf D}_2(4q+2)$, let $a'\in K'$ be as in Lemma 8.6 of \cite{M-R-T}. 
As in the proof of Lemma \ref{Fourth Lemma}, define a homomorphism 
$\eta: \langle C'\rangle \to \langle\ell'\rangle$ as follows. 
Define $\eta(s') = 1$ for all $s'\in C'-\{a'\}$ and define $\eta(a')=\ell'$.  
Then $\eta$ is well defined and $\eta(\ell') = \ell'$. 
By the argument in the proof of Lemma \ref{Fourth Lemma}, the element  
$\ell'$ is in the center of $B^\perp$ and there is a component $L$ of $B^\perp$ 
such that $L\subseteq C$ and $L$ has nontrivial center, 
and if $\ell$ is the longest element of $\langle L\rangle$, 
then $\eta(\ell) = \ell'$ and $\eta(t) = \ell'$ for some $t\in L$. 
Moreover $L$ is of type ${\bf A}_1$, ${\bf C}_{2q+1}$, ${\bf D}_2(4q+2)$, 
${\bf E}_7$, or ${\bf G}_3$ for some $q\geq 1$,  
and if $L$ is of type ${\bf C}_{2q+1}$, then $L\cap{\rm Odd}(t)=\{t\}$. 
As $\eta(z') = 1$, we have that $L\neq K$.  
 
Let $A = \ov{\{t\}}$. 
As $t\in \ov B$, we have that $\ov{\{t\}} \subseteq \ov B$. 
Then $\langle A\rangle$ is conjugate in $\langle C'\rangle$
to $\langle A'\rangle$ for some $A' \subseteq C'$ by Prop. 4.14 of \cite{M-R-T}. 
Now $\eta(\langle A'\rangle) = \eta(\langle A\rangle) = \langle \ell'\rangle$. 
Hence $a'\in A'$. 
Therefore $B'\subseteq A'$ by Lemmas 8.1 and 8.6 of \cite{M-R-T}. 
Hence $B\subseteq A$ by the Basic Matching Theorem, and so $\ov B\subseteq A$. 
Therefore $\ov{\{t\}}=\ov B$. 
As before, $N(t) = B\cup B^\perp$, and so ${\rm Odd}(t) \subseteq L$.  
If $L =\{t\}$, we are done, and  
so we may assume that $L$ is not of type ${\bf A}_1$. 
As before, by applying an automorphism,  
we may assume that $\ell = \ell'$. 

By Lemma 38 of \cite{H-M}, there is an automorphism $\beta$ of $W$ 
such that $\beta(s) = s$ for all $s\in S-(\{r\}\cup {\rm Odd}(t))$, 
and $\beta(r) = r\ell z'$, and $\beta(s) = s\ell z'$ for all $s\in {\rm Odd}(t)$. 
Then $\beta(\ell) = z'$ and $\beta(z') = \ell$. 
As $\beta$ fixes each element of $S-(K\cup L)$ and 
$\beta$ leaves $\langle C\rangle$ invariant, 
we may replace $S$ by $\beta(S)$.  
Then $K$ is replaced by $\beta(L)$,  
and $\beta(K)$ is removed as a possibility for replacing $L$, 
since $\beta(K)$ matches $K'$ and 
$\langle\beta(K)\rangle$ and $\langle K'\rangle$ have the same longest element. 
By the argument in Lemma \ref{Fourth Lemma}, 
we may assume that $\beta(L)$ is of type ${\bf C}_{2q+1}$ or ${\bf D}_2(4q+2)$. 
 
In the above procedure only Coxeter generators of components of $B^\perp$ 
of type ${\bf C}_{2q+1}$, ${\bf D}_2(4q+2)$, ${\bf E}_7$, or ${\bf G}_3$ are replaced. 
By repeating this procedure a finite number of times, 
we can remove the possibility that $L$ is of type 
${\bf C}_{2q+1}$, ${\bf D}_2(4q+2)$, ${\bf E}_7$, or ${\bf G}_3$.  
Then $L=\{t\}$, and we are done.

Assume now that $|\langle B\rangle| > |\langle B'\rangle|$. 
Then $B'$ satisfies conditions (1), (2), (3) of the Blown-Down Theorem. 
Let $r'$ and $K'$ be as in Lemma \ref{Fourth Lemma}. 
If $K'=\{r'\}$, we can blow down $B'$ 
to obtain a set of Coxeter generators $S''$ such that $|S''| = |S'|-1$ 
and $S''$ has one more basic subset than $S'$ isomorphically matching 
basic subsets of $S$, which contradicts the choice of $S'$. 
Therefore $K'$ is of type ${\bf C}_{2q+1}$ or ${\bf D}_2(4q+2)$ for some $q\geq 1$. 
As $N(r') = B'\cup (B')^\perp$, we have that $K'$ is a basic subset of $S'$ 
and $K'\cup (K')^\perp =B'\cup (B')^\perp$.   

Let $K$ be the basic subset of $S$ that matches $K'$. 
Then $K$ isomorphically matches $K'$, since otherwise 
by Theorems 8.4-8.8 of \cite{M-R-T}, we can blow up $S'$ along $K'$   
and then blow down $B'$ 
to obtain a set of Coxeter generators $S''$ such that $|S''| = |S'|$ and $S''$ 
has two more basic subsets than $S'$ isomorphically matching basic subsets of $S$, 
which contradicts the choice of $S'$.

Let $z$ be the longest element of $\langle B\rangle$. 
As in the proof of Lemma \ref{Fourth Lemma}, the element $z$ is in the center of $(B')^\perp$,  
and by applying an automorphism, we may assume 
that $z$ is the longest element of $\langle K'\rangle$. 
Let $\ell$ be the longest element of $\langle K\rangle$. 
If $K$ is of type ${\bf C}_{2q+1}$, let $a\in K$ be such that $K\cap{\rm Odd}(a) = \{a\}$.  
If $K$ is of type ${\bf D}_2(4q+2)$, let $a\in K$ be as in Lemma 8.6 of \cite{M-R-T}. 
As in the proof of Lemma \ref{Fourth Lemma}, define a homomorphism 
$\eta: \langle C\rangle \to \langle\ell\rangle$ as follows. 
Define $\eta(s) = 1$ for all $s\in C-\{a\}$ and define $\eta(a)=\ell$.  
Then $\eta$ is well defined and $\eta(\ell) = \ell$. 
By the argument in the proof of Lemma \ref{Fourth Lemma}, the element  
$\ell$ is in the center of $(B')^\perp$ and there is a component $L'$ of $(B')^\perp$ 
such that $L'\subseteq C'$ and $L'$ has nontrivial center, 
and if $\ell'$ is the longest element of $\langle L'\rangle$, 
then $\eta(\ell') = \ell$ and $\eta(t') = \ell$ for some $t'\in L'$. 
Moreover $L'$ is of type ${\bf A}_1$, ${\bf C}_{2q+1}$, ${\bf D}_2(4q+2)$, 
${\bf E}_7$, or ${\bf G}_3$ for some $q\geq 1$,  
and if $L'$ is of type ${\bf C}_{2q+1}$, then $L'\cap{\rm Odd}(t')=\{t'\}$. 
As $\eta(z) = 1$, we have that $L'\neq K'$.  
 
Let $A' = \ov{\{t'\}}$. 
As $t'\in \ov {B'}$, we have that $\ov{\{t'\}} \subseteq \ov{B'}$. 
Then $\langle A'\rangle$ is conjugate in $\langle C\rangle$
to $\langle A\rangle$ for some $A \subseteq C$ by Prop. 4.14 of \cite{M-R-T}. 
Now $\eta(\langle A\rangle) = \eta(\langle A'\rangle) = \langle \ell\rangle$. 
Hence $a\in A$. 
Therefore $B\subseteq A$ by Lemmas 8.1 and 8.6 of \cite{M-R-T}. 
Hence $B'\subseteq A'$ by the Basic Matching Theorem, and so $\ov{B'}\subseteq A'$. 
Therefore $\ov{\{t'\}}=\ov{B'}$. 
As before, $N(t') = B'\cup (B')^\perp$, and so ${\rm Odd}(t')\subseteq L'$.  
If $L'=\{t'\}$, we derive a contradiction as before.  
Therefore $L'$ is not of type ${\bf A}_1$. 
By applying an automorphism,  
we may assume that $\ell = \ell'$. 

By Lemma 38 of \cite{H-M}, there is an automorphism $\beta$ of $W$ 
such that $\beta(s') = s'$ for all $s'\in S'-(\{r'\}\cup{\rm Odd}(t'))$, 
and $\beta(r') = r'\ell z$, and $\beta(s') = s'\ell z$ for all $s'\in {\rm Odd}(t')$. 
Then $\beta(\ell) = z$ and $\beta(z) = \ell$. 
As $\beta$ fixes each element of $S'-(K'\cup L')$ and 
$\beta$ leaves $\langle C'\rangle$ invariant, 
we may replace $S'$ by $\beta(S')$.  
Then $K'$ is replaced by $\beta(L')$,  
and $\beta(K')$ is removed as a possibility for replacing $L'$, 
since $\beta(K')$ matches $K$ and 
$\langle\beta(K')\rangle$ and $\langle K\rangle$ have the same longest element. 
By the argument in Lemma \ref{Fourth Lemma}, 
we may assume that $\beta(L')$ is of type ${\bf C}_{2q+1}$ or ${\bf D}_2(4q+2)$. 
 
In the above procedure only Coxeter generators of components of $(B')^\perp$ 
of type ${\bf C}_{2q+1}$, ${\bf D}_2(4q+2)$, ${\bf E}_7$, or ${\bf G}_3$ are replaced. 
By repeating this procedure a finite number of times, 
we can remove the possibility that $L'$ is of type 
${\bf C}_{2q+1}$, ${\bf D}_2(4q+2)$, ${\bf E}_7$, or ${\bf G}_3$.  
Then $L'$ is of type ${\bf A}_1$ and we have a contradiction as before. 
Thus the case $|\langle B\rangle| > |\langle B'\rangle|$ leads to a contradiction. 
\end{proof}

\section{The Rank Spectrum of a Coxeter Group}  

In this section, we describe how to determine the set of all possible 
ranks of an arbitrary finitely generated Coxeter group $W$ by inspection 
of any presentation diagram for $W$.

Let $(W,S)$ be a Coxeter system of finite rank. 
Suppose that $S_1,S_2\subseteq S$, with $S=S_1\cup S_2$ and $S_0=S_1\cap S_2$, 
are such that $m(a,b) = \infty$ for all $a\in S_1-S_0$ and $b\in S_2-S_0$. 
Let $\ell\in \langle S_0\rangle$ such that $\ell S_0\ell^{-1} = S_0$.  
The triple $(S_1,\ell,S_2)$ determines an elementary twist of $(W,S)$   
giving a new set of Coxeter generators $S_\ast = S_1\cup\ell S_2\ell^{-1}$ for $W$ 
such that $S_1\cap\ell S_2\ell^{-1} = S_0$. 

Let $B$ be a base of $(W,S)$ of type ${\bf B}_{2p+1}$ or ${\bf D}_2(2p+1)$ for some $p\geq 1$ 
that satisfies the conditions (1), (2), (3) of the Blown-Down Theorem with $K=\{r\}$. 
As $B$ is a simplex, either $B\subseteq S_1$ or $B\subseteq S_2$. 
If $B\subseteq S_1$, define $B_\ast = B$. 
If $B\subseteq S_2$, define $B_\ast =\ell B\ell^{-1}$. 
If $B\subseteq S_0$, then $\ell B\ell^{-1} = B$ by Lemma 4.10 of \cite{M-R-T}, 
and so $B_\ast$ is well defined. 
If $r\in S_1$, define $r_\ast = r$. 
If $r\in S_2$, define $r_\ast = \ell r\ell^{-1}$. 
If $r\in S_0$, then $\ell r\ell^{-1} = r$ by Lemma 4.8 of \cite{M-R-T}, 
and so $r_\ast$ is well defined.

\begin{lemma} 
The set $B_\ast$ is a base of $(W,S_\ast)$ 
of type ${\bf B}_{2p+1}$ or ${\bf D}_2(2p+1)$ for some $p\geq 1$ 
that satisfies the conditions (1), (2), (3) of the Blown-Down Theorem with $\{r_\ast\}$ 
a component of $B_\ast^\perp$. 
\end{lemma}
\begin{proof}
As $B_\ast$ is conjugate to $B$, we deduce that $B_\ast$ is 
a base of $(W,S_\ast)$ of type ${\bf B}_{2p+1}$ or ${\bf D}_2(2p+1)$ for some $p\geq 1$ 
by the Basic Matching Theorem. 
By the Blow-Down Theorem, $W$ has a set of Coxeter generators $S'$ such that $B$ 
matches a base $B'$ of $(W,S')$ with $|\langle B\rangle| < |\langle B'\rangle|$. 
Hence $B_\ast$ matches $B'$ with $|\langle B_\ast\rangle| < |\langle B'\rangle|$. 
Therefore $B_\ast$ satisfies conditions (1) and (2) of the Blow-Down Theorem. 

As $B\cup\{r\}$ is a simplex, either $B\cup\{r\}\subseteq S_1$ or $B\cup\{r\}\subseteq S_2$. 
Hence $r_\ast \in (B_\ast)^\perp$. 
Let $s_\ast\in N(r_\ast)-(B_\ast\cup\{r_\ast\})$. 
Assume first that $B\cup\{r\}\subseteq S_1$. 
Then $B_\ast = B$ and $r_\ast = r$. 
Suppose $s_\ast \in S_1$. 
Then $s_\ast \in N(r) = B\cup B^\perp$. 
Hence $s_\ast \in B^\perp$, and so $s_\ast\in (B_\ast)^\perp$ and $m(r_\ast,s_\ast) = 2$. 
Now suppose that $s_\ast\in \ell S_2\ell^{-1}-S_0$. 
Then $r\in S_0$ and $r=\ell r\ell^{-1}$ and $s_\ast = \ell s\ell^{-1}$ for some $s\in S_2-S_0$. 
Hence $s\in N(r) = B\cup B^\perp$. 
As $B\subset S_1$, we have that $s\in B^\perp$. 
Therefore $m(r,s) = 2$, and so $m(r_\ast,s_\ast) = 2$. 
Moreover, $B\subseteq S_0$, since $s\in S_2-S_0$. 
Hence $B_\ast = \ell B\ell^{-1}$, and so $s_\ast\in (B_\ast)^\perp$. 

Assume now that $B\cup\{r\}\subseteq S_2$. 
Then $B_\ast = \ell B\ell^{-1}$ and $r_\ast = \ell r\ell^{-1}$. 
Suppose $s_\ast \in S_1-S_0$. 
Then $r_\ast \in S_0$, and so $r_\ast = r$. 
Hence $s_\ast \in N(r) = B\cup B^\perp$. 
Now $s_\ast \in B^\perp$, since $s_\ast \in S_1-S_0$. 
Hence $m(r_\ast,s_\ast) = 2$. 
Moreover $B\subseteq S_0$, and so $B_\ast = B$. 
Hence $s_\ast \in (B_\ast)^\perp$. 
Now suppose $s_\ast\in \ell S_2\ell^{-1}$. 
Then $s_\ast = \ell s\ell^{-1}$ for some $s\in S_2$. 
As $1 < m(r_\ast,s_\ast) < \infty$, we have  $1 < m(r,s) < \infty$. 
Hence $s\in N(r) = B\cup B^\perp$. 
Moreover $s\in B^\perp$, since $s_\ast \not\in B_\ast$. 
Hence $m(r,s) = 2$, and so $m(r_\ast,s_\ast) =2$. 
Moreover $s_\ast\in (B_\ast)^\perp$. 
Thus, in all cases, $N(r_\ast) = B_\ast\cup(B_\ast)^\perp$ and $\{r_\ast\}$ 
is a component of $(B_\ast)^\perp$. 
Therefore $B_\ast$ satisfies condition (3) of the Blow-Down Theorem 
with $K = \{r_\ast\}$. 
\end{proof}

Let $B$ be a base of $(W,S)$ of type ${\bf B}_{2p+1}$ or ${\bf D}_2(2p+1)$ for some $p\geq 1$ 
satisfying conditions (1), (2), (3) of The Blow-Down Theorem with $K = \{r\}$. 
We call $r$ a {\it sink} for $B$. 
A element $r$ of $S$ may be a sink for more than one base of $(W,S)$.  
For example, ${\bf A}_1$ is a sink 
for the two bases of ${\bf A}_1\times {\bf D}_2(3)\times {\bf D}_2(3)$. 
A base may have more than one sink. 
For example, the base of ${\bf A}_1\times {\bf A}_1\times {\bf D}_2(3)$ has two sinks.

\begin{lemma} 
Let $B$ be a base of $(W,S)$ of type ${\bf B}_{2p+1}$ or ${\bf D}_2(2p+1)$ for some $p\geq 1$ 
that satisfies the hypothesis of the Blow-Down Lemma with sink $r$. 
Let $S'$ be the set of Coxeter generators obtained by blowing down $S$ along $B$. 
Let $C$ be a base of $(W,S)$ of type ${\bf B}_{2q+1}$ or ${\bf D}_2(2q+1)$ for some $q\geq 1$ 
that satisfies the conditions (1), (2), (3) of the Blow-Down Theorem with sink $s$. 
If $B\neq C$ and $r\neq s$, then $C$ is a base of $(W,S')$ 
that satisfies the conditions (1), (2), (3) of the Blow-Down Theorem with sink $s$. 
\end{lemma}
\begin{proof}
Let $B'$ be the base of $(W,S')$ obtained by blowing down $B$. 
Then $S$ and $S'$ have the same basic subsets except for $B$ and $B'$ 
by Lemma 3.6. 
Therefore $C$ is a base of $(W,S')$. 
By the Blow-Down Theorem, $W$ has a set of Coxeter generators $S''$ 
such that $C$ matches a base $C''$ of $(W,S'')$ with $|\langle C\rangle| < |\langle C''\rangle|$. 
Therefore $C$ satisfies the conditions (1) and (2) of the Blow-Down Theorem. 
As ${\rm Odd}(s) = \{s\}$, we have that $s\not\in B$. 
Therefore $s\in S'$. 

Let $a$ be the element of $S'$ that is not in $S$. 
Then $a = r\ell$ with $\ell$ the longest element of $\langle B\rangle$. 
If $a\not\in N(s)$, then $N(s) = C\cup C^\perp$ 
and $\{s\}$ is a component of $C^\perp$ with respect to $S'$, 
since $s$ is a sink for $C$ with respect to $S$.  
Suppose $a\in N(s)$. 
Then $s\in N(r)$ by Lemma 8.3 of \cite{M-R-T}. 
Hence $s\in B^\perp$. 
Therefore $m(a,s) = 2$. 
As $B\cup\{r\} \subseteq N(s)$, we have that $B\cup\{r\}\subseteq C\cup C^\perp$. 
As $B\neq C$, we have that $B\subseteq C^\perp$.  
As ${\rm Odd}(r) = \{r\}$, we have that $r\not\in C$, and so $r\in C^\perp$. 
Therefore $a \in C^\perp$. 
Hence $N(s) = C\cup C^\perp$ and $\{s\}$ is a component of $C^\perp$ 
with respect to $S'$. 
Thus $C$ satisfies condition (3) of the Blow-Down Theorem, with sink $s$,  
with respect to $S'$. 
\end{proof}

\begin{lemma} 
Let $B$ be a base of $(W,S)$ of type ${\bf C}_{2p+1}$ or ${\bf D}_2(4p+2)$ for some $p\geq 1$ 
along which $(W,S)$ can be blown up.  
Let $S'$ be the set of Coxeter generators obtained by blowing up $S$ along $B$. 
Let $C$ be a base of $(W,S)$ of type ${\bf B}_{2q+1}$ or ${\bf D}_2(2q+1)$ for some $q\geq 1$ 
that satisfies the conditions (1), (2), (3) of the Blow-Down Theorem with sink $s$. 
If $B\neq C$, then $C$ is a base of $(W,S')$ 
that satisfies the conditions (1), (2), (3) of the Blow-Down Theorem with sink $s$. 
\end{lemma}
\begin{proof}
Let $B'$ be the base of $(W,S')$ obtained by blowing up $B$. 
Then $S$ and $S'$ have the same basic subsets except for $B$ and $B'$ 
by Theorems 8.4 and 8.8 of \cite{M-R-T}.  
Therefore $C$ is a base of $(W,S')$. 
By the Blow-Down Theorem, $W$ has a set of Coxeter generators $S''$ 
such that $C$ matches a base $C''$ of $(W,S'')$ with $|\langle C\rangle| < |\langle C''\rangle|$. 
Therefore $C$ satisfies the conditions (1) and (2) of the Blow-Down Theorem. 
As ${\rm Odd}(s) = \{s\}$, we have that $s\not\in B$. 
Therefore $s\in S'$. 

Let $z$ be the longest element of $\langle B\rangle$. 
If $B$ is of type ${\bf C}_{2p+1}$, let $a,b,c$ be the elements of $B$ 
such that $m(a,b) = 4$ and $m(b,c)=3$. Then $N(a)=B\cup B^\perp$ by Theorem 8.2 of \cite{M-R-T}. 
If $B$ is of type ${\bf D}_2(4p+2)$, let $B=\{a,b\}$ with $N(a)=B\cup B^\perp$. 
In either case, let $d =aba$. 
Then $d$ and $z$ are the elements of $S'$ that are not in $S$ 
by Theorems 8.4 and 8.8 of \cite{M-R-T}.

If $d$ and $z$ are not elements of $N(s)$, then $N(s) = C\cup C^\perp$ 
and $\{s\}$ is a component of $C^\perp$ with respect to $S'$, 
since $s$ is a sink for $C$ with respect to $S$.  
Suppose $d$ or $z$ is an element of $N(s)$. 
Then $s\in N(a) = B\cup B^\perp$ by Lemma 8.3 of \cite{M-R-T}. 
Hence $s\in B^\perp$. 
Therefore $m(d,s) = 2$ and $m(z,s)=2$. 
As $B\subseteq N(s)$, we have that $B\subseteq C\cup C^\perp$. 
As $B\neq C$, we have that $B\subseteq C^\perp$. 
Therefore $\{d,z\}\subseteq C^\perp$. 
Hence $N(s) = C\cup C^\perp$ and $\{s\}$ is a component of $C^\perp$ 
with respect to $S'$. 
Thus $C$ satisfies condition (3) of the Blow-Down Theorem, with sink $s$,  
with respect to $S'$. 
\end{proof}

\begin{lemma} 
Let $B$ be a base of $(W,S)$ of type ${\bf C}_{2p+1}$ or ${\bf D}_2(4p+2)$ for some $p\geq 1$ 
along which $(W,S)$ can be blown up.  
Let $S'$ be the set of Coxeter generators obtained by blowing up $S$ along $B$. 
Let $C$ be a base of $(W,S)$ of type ${\bf C}_{2q+1}$ or ${\bf D}_2(4q+2)$ for some $q\geq 1$ 
along which $(W,S)$ can be blown up. 
If $B\neq C$, then $C$ is a base of $(W,S')$ 
along which $(W,S')$ can be blown up. 
\end{lemma}
\begin{proof}
Let $B'$ be the base of $(W,S')$ obtained by blowing up $B$. 
Then $S$ and $S'$ have the same basic subsets except for $B$ and $B'$ 
by Theorems 8.4 and 8.8 of \cite{M-R-T}. 
Therefore $C$ is a base of $(W,S')$.

Let $z$ be the longest element of $\langle B\rangle$. 
If $B$ is of type ${\bf C}_{2p+1}$, let $a,b,c$ be the elements of $B$ 
such that $m(a,b) = 4$ and $m(b,c)=3$. Then $N(a)=B\cup B^\perp$ by Theorem 8.2 of \cite{M-R-T}. 
If $B$ is of type ${\bf D}_2(4p+2)$, let $B=\{a,b\}$ with $N(a)=B\cup B^\perp$. 
In either case, let $d =aba$. 
Then $d$ and $z$ are the elements of $S'$ that are not in $S$ 
by Theorems 8.4 and 8.8 of \cite{M-R-T}.  

Let $v$ be the element of $C$ such that $N(v) = C\cup C^\perp$ 
as in Theorems 8.5 and 8.7 of \cite{M-R-T}. 
If $d$ and $z$ are not elements of $N(v)$, then $N(v) = C\cup C^\perp$ with respect to $S'$.  
Suppose $d$ or $z$ is an element of $N(v)$. 
Then $v\in N(a) = B\cup B^\perp$ by Lemma 8.3 of \cite{M-R-T}. 
Hence $v\in B^\perp$. 
As $B\subseteq N(v)$, we have that $B\subseteq C\cup C^\perp$. 
As $B\neq C$, we have that $B\subseteq C^\perp$. 
Therefore $\{d,z\}\subseteq C^\perp$. 
Hence $N(v) = C\cup C^\perp$ with respect to $S'$. 
Then $(W,S')$ can be blown up along $C$ by Theorems 8.4 and 8.8 in \cite{M-R-T}. 
\end{proof}

\begin{theorem} {\rm (Rank Spectrum Theorem)} 
Let $\{B_1,\ldots, B_k\}$, $k\geq 0$, be a maximal set of bases of $(W,S)$ 
of type ${\bf B}_{2p+1}$ or ${\bf D}_2(2p+1)$ for some $p\geq 1$ 
that satisfy the conditions of the Blown-Down Theorem  
with distinct sinks $\{s_1,\ldots, s_k\}$. 
Let $C_1,\ldots, C_\ell$, $\ell \geq 0$,
be the bases of $(W,S)$ of type ${\bf C}_{2q+1}$ or ${\bf D}_2(4q+2)$ 
for some $q\geq 1$ along which $(W,S)$ can be blown up. 
Then the set of all possible ranks of $W$ is $\big\{|S|-k, \ldots, |S|+\ell\big\}$. 
\end{theorem}
\begin{proof}
By Lemmas 5.1 and 5.2, we get a sequence of sets of Coxeter generators $S_0,\ldots, S_k$ for $W$ 
such that $S = S_0$ and $S_i$ is obtained from $S_{i-1}$ by twisting $S_i$, as in Theorem 3.7,  
and then blowing down along a base conjugate to $B_i$ for each $i=1,\ldots,k$. 
Then $S_k$ has minimum rank over all sets of Coxeter generators of $W$ by Theorem 4.1 and Lemma 5.3. 
Hence, the minimum rank of $W$ is $|S|-k$.  

Let $a_i \in C_i$ be the element of $C_i$, for $i=1,\ldots,\ell$ 
that is removed in the blowing up process. 
As $N(a_i) = C_i\cup C_i^\perp$ for each $i$, we have that $C_i$ 
is the only base of $(W,S)$ that contains $a_i$ for each $i$. 
By Lemma 5.4, we have a sequence of sets of 
Coxeter generators $S^{(0)},\ldots, S^{(\ell)}$ for $W$ such that 
$S = S^{(0)}$ and $S^{(i)}$ is obtained from $S^{(i-1)}$ by blowing up $S^{(i)}$ along $C_i$ 
for each $i=1,\ldots,\ell$. 
Then $S^{(\ell)}$ has maximum rank over all sets of Coxeter generators of $W$ 
by Theorem 9.1 of \cite{M-R-T}. 
Hence, the maximum rank of $W$ is $|S|+\ell$. 
Thus the set of all possible ranks of $W$ is $\big\{|S|-k, \ldots, |S|+\ell\big\}$. 
\end{proof}

The numbers $k$ and $\ell$ in the Rank Spectrum Theorem 
can be determined by inspecting the presentation diagram of $(W,S)$. 
For example, $k =\ell = 1$ for the system 
${\bf A}_1\times {\bf D}_2(3)\times{\bf D}_2(6)$.

\end{document}